\documentclass[12pt]{article}
\usepackage{amssymb}
\usepackage{amsmath}
\usepackage{amsthm}
\usepackage{color}
\usepackage{comment}
\usepackage{geometry}		
\usepackage{graphicx}
\usepackage{float}
\usepackage{booktabs,longtable}
\usepackage[english]{babel}
\usepackage{multirow}
\usepackage{tabularx,booktabs}
\usepackage{caption}
\usepackage{subcaption}
\usepackage{hyperref}
\hypersetup{
    colorlinks=true,
    linkcolor=blue,
    urlcolor=blue,
    citecolor=blue
}

\usepackage{mathtools}
\usepackage{tensor}

\theoremstyle{definition}

\theoremstyle{remark}

\usepackage{authblk} 
\usepackage{orcidlink} 
\usepackage{hyperref}
\hypersetup{colorlinks=true, linkcolor=blue, urlcolor=blue}

\usepackage{multicol}
\usepackage{changepage}
\let\originalleft\left
\let\originalright\right
\renewcommand{\left}{\mathopen{}\mathclose\bgroup\originalleft}
\renewcommand{\right}{\aftergroup\egroup\originalright}

\usepackage{tabularx}
\usepackage{caption}
\usepackage{subcaption}

\usepackage{tikz}
\usetikzlibrary{arrows.meta, positioning,backgrounds}

\title{Stochastic Modelling and Analysis of Within-Farm Highly Pathogenic Avian Influenza Dynamics in Dairy Cattle}

\author[1]{P.~Tiwari\,\orcidlink{0000-0003-0206-6604}\thanks{Corresponding authors: \href{mailto:parul.tiwari@aut.ac.nz}{parul.tiwari@aut.ac.nz},~\href{mailto:hammed.fatoyinbo@aut.ac.nz}{hammed.fatoyinbo@aut.ac.nz}}}
\author[1]{M.~Smitha}
\author[1]{H.O.~Fatoyinbo\,\orcidlink{0000-0002-6036-2957}*}

\affil[1]{Department of Mathematical Sciences, Auckland University of Technology, Auckland 1010, New Zealand}
\begin{document}
\maketitle


\begin{abstract}
Highly pathogenic avian influenza (HPAI) has expanded its host range with recent detections in dairy cattle, raising critical concerns regarding within-herd persistence and cross-species spillover. This study develops a stochastic $SEI_sI_aR-B$ compartmental model to analyse HPAI transmission, explicitly accounting for environmental pathogen reservoirs and noise intensities through Wiener processes. The positivity and boundedness of solutions are established, and the disease-free and endemic equilibria are analytically derived. The basic reproduction number is determined using the next-generation matrix method. Numerical simulations confirm that the model dynamics are consistent with theoretical analysis and illustrate how stochastic fluctuations significantly influence disease persistence. Furthermore, sensitivity analysis using Latin Hypercube Sampling (LHS) and Partial Rank Correlation Coefficients (PRCC) identifies the transmission rate from asymptomatic infectious cattle ($\beta_a$) as the primary driver of transmission. The model effectively captures the dynamics of environmental variability affecting HPAI spread, suggesting that effective control strategies must prioritise the early detection and isolation of asymptomatic carriers alongside environmental management.
\end{abstract}

\section{Introduction}

Avian influenza (AI), commonly referred to as bird flu, is a highly contagious viral disease caused by influenza A viruses belonging to the family Orthomyxoviridae and the genus Influenza A \cite{opata2025predictiveness,yamaji2020pandemic}. These viruses primarily infect avian species, especially wild aquatic birds, which act as the natural reservoir and contribute to global transmission through migratory movements \cite{Mao2024, FaPRVic25}. While wild birds are the primary host, spillover into domestic poultry can lead to widespread outbreaks, causing severe economic and ecological consequences. Avian influenza viruses are classified into two categories: low pathogenic avian influenza (LPAI), which causes mild or subclinical infections, and highly pathogenic avian influenza (HPAI), which results in severe systemic disease and high mortality rates \cite{Rodriguez2024JDSC}.

Globally, HPAI H5N1 subtype remains one of the most concerning strains due to its ability to infect multiple species and its persistence in the environment. Currently, only New Zealand, Australia, and Antarctica remain free from HPAI H5N1 \cite{docnz, stanislawek2024}, outbreaks in other parts of the world highlight the virus’s capacity for rapid adaptation and cross-species transmission. This cross-species potential has led to heightened global vigilance in monitoring wild bird populations, their environment and domestic animal interfaces.

The recent emergence of HPAI H5N1 in non-avian hosts, particularly in dairy cattle, marks a significant epidemiological shift \cite{BUTT2024S13}. In 2024, the United States confirmed several outbreaks of H5N1 in dairy herds, with detections in nasal swabs and unpasteurized milk, suggesting viral replication within the mammary gland and possible environmental contamination through the milking system \cite{APHIS2024Brief,Caserta2024Nature,usda2024,cdc2024,Rodriguez2024JDSC}. The initial cases were detected in Texas and subsequently in Kansas, Michigan, and New Mexico. These findings challenge the long-held assumption that cattle are resistant to avian influenza infection and highlight new pathways for viral maintenance outside avian populations \cite{Kaiser2025EID, Rawson2025}. The transmission was linked to direct cow-to-cow contact, contaminated milking equipment, human-assisted spread, and potential exposure to infected wild birds or poultry. Clinically, affected cattle exhibited reduced milk production, fever, nasal discharge, and, in many cases, mild or asymptomatic presentations, complicating detection \cite{mbio2024,PenaMosca2025}.

The incubation phase of HPAI in cattle can last up to two weeks; during this window, the movement of asymptomatic carriers between farms creates a risk of further geographic dissemination \cite{incubation}. Although H5N1 is highly lethal in avian species, mortality and culling rates in cattle have historically remained low, typically staying below 2\% \cite{avma2024}. While most infected cattle recover with supportive treatment, the virus continues to present substantial economic and biosecurity challenges \cite{Roche2020JDS}. Immunity gained from prior infection is expected to be partial or transient—particularly given viral evolution—making lifelong protection improbable \cite{Facciuolo2025}. Consequently, without rigorous surveillance and biosecurity protocols, recurrent outbreaks remain a distinct possibility \cite{Mrope2024}.

Avian influenza transmission occurs through various pathways, including direct interaction with infected hosts, inhalation or ingestion of viral particulates, contact with contaminated fomites, and exposure to environmental reservoirs \cite{transmission,Mrope2024}. While bird-to-bird transmission is well documented, the identification of HPAI in mammals raises concerns regarding environmental persistence and potential adaptation to mammalian physiology. Occupational groups such as veterinarians, poultry and dairy workers, and animal handlers face the highest risk of exposure. In birds, clinical signs often include neurological deficits—such as tremors, torticollis, and ataxia—alongside respiratory distress and diarrhoea, with severe cases frequently resulting in sudden death. Although human infections are rare, they can manifest as fever, muscle pain, and respiratory issues, potentially progressing to pneumonia or multi-organ failure \cite{PenaMosca2025}.

The emergence of HPAI H5N1 in cattle represents a major shift in the virus’s epidemiological behavior, with significant consequences for animal health, public health, and the agricultural industry \cite{Rodriguez2024JDSC,APHIS2024Brief}. Effective control strategies require a detailed understanding of how the virus spreads within herds--especially the relative contribution of direct animal-to-animal transmission compared with environmental pathways, and the conditions that enable outbreaks to continue \cite{Mrope2024,Regassa2024,math13193086}. Mathematical modelling offers a powerful framework for investigating these issues by quantifying distinct transmission modes, pinpointing key epidemiological thresholds, and assessing the likely impact of different control options \cite{Wu2013Sensitivity,Heffernan2005, FaTiPGh25}. As the epidemiology of HPAI evolves, there is a pressing need for mathematical models that explicitly incorporate both direct and environmental transmission processes. These models are essential for estimating infection burdens, determining critical control levers, and comparing the performance of intervention strategies across a range of epidemiological settings. Recently, the transmission dynamics of Avian influenza in cattle population in the United States has been studied using mathematical models in \cite{Rawson2025}.

Motivated by the recent HPAI H5N1 outbreak in dairy cattle in the US \cite{cdc2024, Rawson2025}, we develop and analyse a stochastic within-farm model for HPAI transmission. The model incorporates both symptomatic and asymptomatic cattle, as well as virus contamination in the environment. Using standard stochastic differential equation theory, we establish the positivity and boundedness of solutions and derive conditions for the extinction and persistence of HPAI. Numerical experiments are carried out to support and illustrate the theoretical findings.

The remainder of the paper is structured as follows: Section~\ref{sec:modelformulation} presents the formulation of the stochastic HPAI model. Section~\ref{sec:modelanalysis} is devoted to the qualitative analysis of the model. Numerical simulations are reported in Section~\ref{sec:numel}. Finally, the conclusion and study limitations are given in Sections~\ref{sec:conclusion} and~\ref{sec:limitation}, respectively.



\section{Model formulation}\label{sec:modelformulation}

The following key assumptions guide the model formulation:
\begin{itemize}
    \item The cattle herd is treated as an open population: animals enter through recruitment, purchase, or births at a constant rate $\Lambda$, and leave through natural mortality at rate $\mu$ or disease-induced mortality at rate $d$. This setup reflects realistic herd management practices rather than a fully closed herd.

\item Two infection pathways are considered: (a) direct transmission, when susceptible cattle have contact with infectious animals at the rate $\beta_a$ and $\beta_{s}$, and (b) indirect transmission at rate $\beta_{B}$, when cattle encounter virus particles present in the environment (such as on contaminated surfaces, bedding, or milking equipment).

\item Infected animals may develop clinical signs or remain asymptomatic. Both groups shed virus, but asymptomatic cattle typically go unnoticed, prolonging their infectious period and helping to maintain transmission within the herd. 

\item Once infected animals are identified and removed from the herd through culling, separation, slaughter, or other means, they are assumed no longer to contribute to transmission and are placed in the removed class ($R$).

\item The environmental compartment ($B$) represents viral particles deposited into the environment by infected cattle. This environmental burden declines at the constant rate $\varepsilon$, reflecting loss of infectivity or degradation under natural conditions (e.g., temperature, sunlight, disinfectants). 

\item Infection arising from environmental exposure is modelled using a saturating function, so that the risk of infection increases with contamination but levels off once exposure becomes very high.
\end{itemize}


Based on the assumptions, a mathematical model named $SEI_{s}I_{a}R-B$ for HPAI was developed where total cattle population is divided into five compartments: Susceptible ($S$), Exposed ($E$), Symptomatic ($I_{s}$), Asymptomatic ($I_{a}$), Removed ($R$), and contaminated environment ($B$). The model equations are written using a system of ordinary differential equations (ODEs) and given by
\begin{equation}
\begin{aligned}
\frac{dS}{dt} &= \Lambda-\beta_s S\frac{I_s}{N} -\beta_a S\frac{I_a}{N} -\beta_B S\frac{B}{K+B}-\mu S, \\
\frac{dE}{dt} &= \beta_sS\frac{I_s}{N}+ \beta_a S\frac{I_a}{N} +\beta_B  S\frac{B}{K+B}-\sigma E-\mu E, \\ 
\frac{dI_s}{dt} &= \nu\sigma E-(\mu+d+\gamma)I_s, \\ 
\frac{dI_a}{dt} &= (1-\nu)\sigma E-(\mu + d+\delta)I_a,\\ 
\frac{dR}{dt} &= \gamma I_s+\delta I_a-\mu R,\\
\frac{dB}{dt} &= \omega_s I_s+\omega_a I_a-\varepsilon B.
\label{eq:ode-full}
\end{aligned}
\end{equation}
with initial conditions
$$S(0)\ge 0, E(0)\ge 0, I_{s}(0)\ge 0, I_{a}(0)\ge 0, R(0)\ge 0, B(0)\ge 0.$$
The total population at any time $t$ satisfies $$N(t) = S + E + I_s + I_a + R.$$

\begin{figure}[h!]
    \centering
    \begin{tikzpicture}[scale=0.7, transform shape, node distance=2.8cm, auto, >=Stealth, font=\large]
    \node [draw, circle, minimum size=1cm, line width= 2pt, fill=blue!50] (S) {$S$};
     \node [draw, circle, minimum size=1cm, line width= 2pt, right=of S,fill=yellow!80] (E) {$E$};
     \node[draw, circle, minimum size=1cm, line width= 2pt, above right=1.5cm and 2cm of E, fill=orange!80] (Ia) {$I_a$};
     \node[draw, circle, minimum size=1cm, line width= 2pt, below right=1.5cm and 2cm of E, fill=red!80] (Is) {$I_s$};
     \node [draw, circle, minimum size=1cm, line width= 2pt, right=4cm of E,fill=green!80] (R) {$R$};

     \node [draw, circle, minimum size=1cm, line width= 2pt, below=of S,xshift=25mm, fill=magenta!80] (B) {$B$};

    \draw [->, line width = 0.7mm] (S) -- (E) node[midway, above] {$\lambda$};
    \draw[->, line width = 0.7mm] (E) -- node[above left] {$\sigma(1-\nu)$} (Ia);
 \draw[->, line width = 0.7mm] (E) -- node[above right] {$\sigma \nu$} (Is);
 \draw[->, line width = 0.7mm] (Ia) -- node[above right] {$\delta$} (R);
 \draw[->, line width = 0.7mm] (Is) -- node[below right] {$\gamma$} (R);

\draw[->, dashed, line width=0.6pt] 
    (Is) -- node[midway, above] {$\omega_s$} (B);

\begin{pgfonlayer}{background} 
  \draw[->, dashed, line width=0.6pt]
    (Ia) .. controls + (3.5,1.8) and + (15.5,-5.2) .. (B) node[midway, right] {$\omega_a$};
\end{pgfonlayer}

\path (S) -- (E) coordinate[midway] (midSE);
\draw[->, dashed, line width=0.6pt] 
    (B) -- node[midway, left] {$\beta_B$} (midSE);
    
     \path[->, line width = 0.7mm]  (S) edge node[swap, yshift=-7mm,xshift=3mm] {$ \mu$} +(0,-1.6);
     \path[->, line width = 0.7mm]  (E) edge node[swap, yshift=-7mm,xshift=3mm] {$\mu$} +(0,-1.5);
     \path[->, line width = 0.7mm] (Is) edge node[swap, yshift=-7mm,xshift=9mm] {$(\mu+d)$} +(0,-1.5);
     \path[->, line width = 0.7mm] (Ia) edge node[swap, yshift=7mm, xshift=-10mm] {$(\mu+d)$} +(0,1.5);;
     \path[->, line width = 0.7mm]  (R) edge node[swap, yshift=-7mm,xshift=3mm] {$ \mu$} +(0,-1.6);
       \path[->, line width = 0.7mm]  (B) edge node[swap, yshift=-7mm,xshift=3mm] {$\varepsilon$} +(0,-1.5);
  
 \draw[->, line width = 0.7mm] (-2,0) node[left] {$\Lambda$} -- (S);
  \end{tikzpicture}
 \caption{Flow diagram illustrating the transmission routes of the HPAI virus in cattle}
 \label{fig:SEIASRB}
\end{figure}
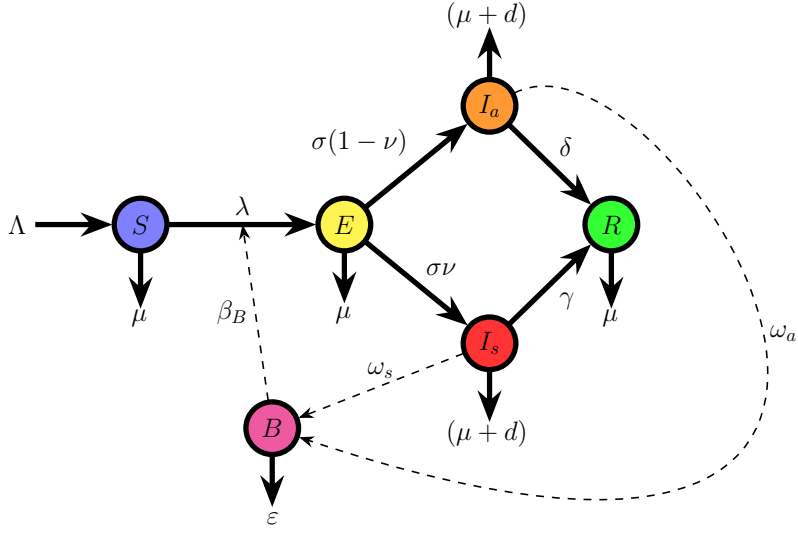

A schematic diagram of the model structure is presented in Figure~\ref{fig:SEIASRB}, illustrating the flow of individuals between compartments and the interactions between cattle and their environment. Each arrow denotes a transition governed by a specific rate parameter, such as infection, progression, recovery, or decay.

To capture the inherent randomness and environmental variability influencing the transmission dynamics of avian influenza (AI), we extend the deterministic compartmental framework to a stochastic differential equation (SDE) system by incorporating multiplicative white noise perturbations in each epidemiological class except recovery class as it will not affect the dynamics of infected class \cite{math11194199}. The deterministic model assumes perfect homogeneity in disease transmission and parameter constancy over time; however, in natural ecological settings such as cattle migration intensity, pathogen survival in the environment, and host immune variability can significantly alter infection rates and removal processes \cite{Allen2017IDM, gray2011stochastic}. The stochastic version of the model (\ref{eq:sde-full}) therefore provides a more realistic representation of AI dynamics by allowing such uncertainties to enter the system as diffusion terms proportional to the state variables.

\begin{equation}
\begin{aligned}
\frac{dS}{dt} &= \Lambda-\beta_s S\frac{I_s}{N} -\beta_a S\frac{I_a}{N} -\beta_B S\frac{B}{K+B}-\mu S+ \sigma_S S\, \frac{dW_S}{dt}\\ \frac{dE}{dt} &= \beta_sS\frac{I_s}{N}+ \beta_a S\frac{I_a}{N} +\beta_B  S\frac{B}{K+B}-\sigma E-\mu E + \sigma_e E\, \frac{dW_e}{dt}\\ \frac{dI_s}{dt} &= \nu\sigma E-(\mu+d+\gamma)I_s + \sigma_s I_s\, \frac{dW_s}{dt}\\ \frac{dI_a}{dt} &= (1-\nu)\sigma E-(\mu + d+\delta)I_a + \sigma_a I_a\, \frac{dW_a}{dt}\\ \frac{dR}{dt} &= \gamma I_s+\delta I_a-\mu R\\ \frac{dB}{dt} &= \omega_s I_s+\omega_a I_a-\varepsilon B + \sigma_B B\, \frac{dW_B}{dt} 
\label{eq:sde-full}
\end{aligned}
\end{equation}

Each compartment is influenced by an independent Wiener process $W_i(t),\, i = S, e, s, a, B$ representing the cumulative effect of random perturbations acting on that class. The term $\sigma_S S\, dW_S(t)$ captures random fluctuations in the recruitment and natural death processes of susceptible hosts due to unpredictable environmental events that changes host availability or contact rates. The exposed class includes the diffusion term $\sigma_e E\, dW_e(t)$, which reflects uncertainty in the latent period, for instance, due to temperature-dependent virus incubation to infected excreta. The stochastic perturbation $\sigma_s I_s\, dW_s(t)$ in the symptomatic infected class models random variability in transmission and recovery among visibly infected individuals, possibly caused by variations in immune response. Similarly, $\sigma_a I_a\, dW_a(t)$ in the asymptomatic infected compartment accounts for random changes in subclinical infection dynamics, such as variations in viral shedding or host contact behavior that are difficult to monitor in the field. The environmental virus reservoir $B(t)$ is perturbed by $\sigma_B B\, dW_B(t)$, representing unpredictable changes in viral persistence and decay in the environment, affected by climatic and ecological factors \cite{chang2025long}.

Each noise term is multiplicative which shows that the magnitude of random fluctuations scales with the population level of the corresponding compartment. This is biologically feasible because the stochastic influence vanishes when the compartment population is zero and increases with population density. Such a formulation guarantees positivity of solutions and preserves the nonnegativity of all compartments \cite{mao2007stochastic}. The stochastic framework thus enables investigation of how random environmental shocks and intrinsic variability can drive disease persistence, extinction, or oscillations even when deterministic conditions predict eradication. The parameters used in the model and their biological interpretation are summarised in Table~\ref{tab:parameters}. 

\begin{table}[H]
\centering
\caption{Model parameters and their biological interpretation}
\begin{tabular}{|c | p{10cm}|}
\hline
\textbf{Parameter} & \textbf{Definition} \\
\hline
$\Lambda$ & Recruitment rate of susceptible cattle into the population \\
$\mu$ & Natural mortality rate of cattle \\
$\beta_s$ & Transmission rate from symptomatic infectious cattle to susceptible cattle via direct contact \\
$\beta_a$ & Transmission rate from asymptomatic infectious cattle to susceptible cattle via direct contact \\
$\beta_B$ & Transmission rate from the contaminated environment to susceptible cattle \\
$K$ & Half-saturation constant representing the environmental pathogen load at which transmission rate is half-maximal \\
$\sigma$ & Rate at which exposed cattle progress to the infectious stage \\
$\nu$ & Proportion of exposed cattle that become symptomatic upon infection \\
$\gamma$ & Removal rate of symptomatic cattle (due to recovery, isolation, culling, or slaughter) \\
$\delta$ & Removal rate of asymptomatic cattle (due to recovery, isolation, culling, or slaughter) \\
$d$ & Disease-induced mortality rate in symptomatic cattle \\
$\omega_s$ & Rate of virus shedding from symptomatic cattle into the environment \\
$\omega_a$ & Rate of virus shedding from asymptomatic cattle into the environment \\
$\varepsilon$ & Rate of environmental virus decay or removal \\
\hline
\end{tabular}
\label{tab:parameters}
\end{table}

\section{Analysis of model} \label{sec:modelanalysis}
\subsection{Positivity and Boundedness}
Consider the stochastic AI system in Itô form
\begin{equation}
\left.
\label{eq:stochastic-system}
\begin{aligned}
	dS &= \Big[\Lambda-\beta_s\tfrac{SI_s}{N}-\beta_a\tfrac{SI_a}{N}
	-\beta_B\tfrac{SB}{K+B}-\mu S\Big]\,dt + \sigma_S S\, dW_S(t)\\
	dE &= \Big[\beta_s\tfrac{SI_s}{N}+\beta_a\tfrac{SI_a}{N}
	+\beta_B\tfrac{SB}{K+B}-\sigma E-\mu E\Big]\,dt + \sigma_E E\, dW_e(t)\\
	dI_s &= \big[\nu\sigma E-(\mu+d+\gamma)I_s\big]\,dt + \sigma_s I_s\, dW_s(t)\\
	dI_a &= \big[(1-\nu)\sigma E-(\delta+\mu+d)I_a\big]\,dt + \sigma_a I_a\, dW_a(t)\\
	dR &= (\gamma I_s+\delta I_a-\mu R)\,dt\\
	dB &= (\omega_s I_s+\omega_a I_a-\varepsilon B)\,dt + \sigma_B B\, dW_B(t)
\end{aligned}
\right\}
\end{equation}
where $W_i(t)$ are independent standard Brownian motions.
\subsubsection{Positivity}
Assuming the standard hypotheses that all parameters are nonnegative and the initial state
	\[
	(S_0,E_0,I_{s,0},I_{a,0},R_0,B_0)\in\mathbb{R}_+^6
	\]
	is nonnegative a.s. Then the Itô SDE system admits a unique global strong solution \, \, \(X_t=(S_t,E_t,I_{s,t},I_{a,t},R_t,B_t)\) with
	\[
	X_t \in \mathbb{R}_+^6 \quad\text{for all }t\ge0\quad\text{a.s.,}
	\]
	and the first moments satisfy uniform-in-time bounds and hence the solution is bounded.
\begin{proof}
Since the diffusion terms are proportional to the state variables 
($\sigma_i X_i\, dW_i$), we can apply Itô’s formula to each class of model \eqref{eq:stochastic-system} for every strictly positive component $X_t$. Writing SDE in compact form as
\begin{equation}
   dX_t = f_X(X_t)\,dt + \sigma_X X_t\,dW_t
   \label{eq:Ito SDE}
\end{equation}
    where \(f_X(\cdot)\) denotes the drift of that component and is locally Lipschitz.
	Apply Itô’s formula to \(\ln X_t\) while \(X_t>0\):
\begin{equation}
	d(\ln X_t) = \left(\frac{f_X(X_t)}{X_t} - \tfrac{1}{2}\sigma_X^2\right)\,dt + \sigma_X\,dW_t.
\label{eq:logIto SDE}
\end{equation}
    The right-hand side of \eqref{eq:logIto SDE} is finite for all finite times before explosion because \(f_X\) is locally bounded on compacts and \(X_t>0\) on that interval.

    This means, for any bounded subset $P \subset \mathbb{R}_+^6$,
    there exists a constant $C_P>0$ such that
    \[
    |f_X(x)| \le C_P, \qquad \forall\, x\in P.
    \]
    Thus, whenever the state variables
    $(S_t,E_t,I_{s,t},I_{a,t},R_t,B_t)$ remain within a finite range,
    the drift term $f_X$ cannot become infinite.
    Since all drift terms in the stochastic AI model are composed of
    polynomial and rational expressions with positive denominators,
    they are locally bounded on the positive orthant $\mathbb{R}_+^6$.
    This ensures that the stochastic integrals
    $\int_0^t f_X(X_z)\,dz$ are well defined and finite for all finite $t$.

Now, When applying Itô’s formula to $\ln X_t$,
    it is necessary that $X_t>0$, so that $\ln (X_t)$ is well defined.
    Therefore, the analysis is carried out on a maximal time interval
    $[0,\tau)$ for which $X_t>0$ almost surely.
    If we assume, for contradiction, that $X_t$ reaches zero at a finite time $\tau$,
    then $\ln (X_t) \to -\infty$ as $t\to\tau$.
    However, since $f_X$ is locally bounded and $\sigma_X$ is finite,
    the right–hand side of the Itô differential \eqref{eq:logIto SDE} remains finite on $[0,\tau)$, leading to a contradiction.
    Hence $X_t$ cannot reach zero in finite time and remains positive almost surely.
	
	Also, the recovered class \(R_t\) has no diffusion term and its drift is \(\gamma I_{s,t}+\delta I_{a,t}-\mu R_t\). Since \(I_{s,t},I_{a,t}\ge0\) a.s., standard ODE comparison shows \(R_t\) cannot become negative even if \(R_0\ge0\).
	
	Hence all components remain strictly positive whenever their initial value is positive for all \(t\ge0\) almost surely.
\end{proof}
\subsubsection{Boundedness}
We have
	\[
	N_t := S_t + E_t + I_{s,t} + I_{a,t} + R_t.
	\]
	Summing the SDEs \eqref{eq:stochastic-system} for \(S_t,E_t,I_{s,t},I_{a,t},R_t\) yields
	\[
	dN_t = \big(\Lambda - \mu N_t - d I_{s,t}\big)\,dt
	+ \sigma_S S_t\,dW_S + \sigma_e E_t\,dW_e + \sigma_s I_{s,t}\,dW_s + \sigma_a I_{a,t}\,dW_a.
	\]
	Taking expectation and using \(\mathbb{E}[dW_i]=0\), we get
	\[
	\frac{d}{dt}\mathbb{E}[N_t] = \Lambda - \mu \mathbb{E}[N_t] - d\,\mathbb{E}[I_{s,t}]
	\le \Lambda - \mu \mathbb{E}[N_t]
	\]
	Solving this linear differential inequality, we get
	\[
	\mathbb{E}[N_t] \le \mathbb{E}[N_0] e^{-\mu t} + \frac{\Lambda}{\mu}\big(1-e^{-\mu t}\big)
	\le \max\{\mathbb{E}[N_0],\Lambda/\mu\},
	\]
	so the host total is uniformly bounded in expectation.
	
	Now, we consider the environmental compartment \(B_t\). Using SDE from \eqref{eq:stochastic-system},
	\[
	\frac{d}{dt}\mathbb{E}[B_t] = \mathbb{E}[\omega_s I_{s,t} + \omega_a I_{a,t} - \varepsilon B_t]
	\]
	Using the bound \(I_{s,t}+I_{a,t}\le N_t\) and let \(\omega_{\max}=\max\{\omega_s,\omega_a\}\),
	\[
	\frac{d}{dt}\mathbb{E}[B_t] \le \omega_{\max}\,\mathbb{E}[N_t] - \varepsilon \mathbb{E}[B_t]
	\]
	Since \(\mathbb{E}[N_t]\) is uniformly bounded, this is a linear inhomogeneous differential inequality for \(\mathbb{E}[B_t]\). Solving using variation of constants gives
	\[
	\mathbb{E}[B_t] \le \mathbb{E}[B_0] e^{-\varepsilon t} + \frac{\omega_{\max}\sup_{s\ge0}\mathbb{E}[N_s]}{\varepsilon}\big(1-e^{-\varepsilon t}\big),
	\]
	hence \(\sup_{t\ge0}\mathbb{E}[B_t]<\infty\)
	
	Thus, every component of model \eqref{eq:stochastic-system} has a finite first moment uniformly in time. Using standard SDE theory \cite{khasminskii2012stochastic}, bounded first moments prevent finite-time explosion and thus the local solution extends to all \(t\ge0\).
    
\subsection{Equilibrium analysis}
\subsubsection{Disease-free equilibrium (DFE)}
At the disease-free equilibrium (DFE), the infection is completely absent from the cattle population, meaning that no individuals exist in the exposed or infectious compartments. To obtain this equilibrium, we set all time derivatives in the system to zero and assume that no disease transmission occurs ($E = I_s = I_a = R = B = 0$). Under these conditions, the number of susceptible cattle remains constant, determined only by the balance between recruitment and natural death. The noise term in stochastc model does not contribute to equilibrium itself because its expectation is zero ($\mathop{\mathbb{E}(dW_i(t) =0}$).
However, noise does affect the stability and long-term behaviour around that equilibrium \cite{chang2025long}.
At the disease–free equilibrium all infected compartments vanish,\\
so set \(E^*=I_s^*=I_a^*=B^*=0\). Thus we have
\[
0 = \Lambda - \mu S^*,\qquad 0 = -\mu R^*.
\]
Hence the DFE is
\[
 \; \text{DFE} \;=\; \Big(S^*,E^*,I_s^*,I_a^*,R^*,B^*\Big)
\;=\; \Big(\tfrac{\Lambda}{\mu},\;0,\;0,\;0,\;0,\;0\Big). \;
\]

\subsubsection{Basic reproduction number \(\mathcal{R}_0\)}
The basic reproduction number, denoted as $\mathcal{R}_0$, represents the average number of secondary infections produced by a single infectious cattle introduced into a completely susceptible population \cite{Heffernan2005}. It serves as a key threshold parameter that determines whether an infection can invade and persist within a population.

\noindent We derive \(\mathcal{R}_0\) using the next–generation matrix denoted by $\psi_{NGM}$. At the DFE we have \(N^*=S^*\), hence \(\dfrac{S^*}{N^*}=1\). Therefore the linearised new–infection matrix \(F\) at DFE is
\[
F \;=\;
\begin{pmatrix}
0 & \beta_s & \beta_a & \dfrac{\beta_B S^*}{K} \\
0 & 0      & 0      & 0\\
0 & 0      & 0      & 0\\
0 & 0      & 0      & 0
\end{pmatrix}.
\]

\noindent The transition matrix \(V\) describes the flow of individuals out of infected compartments due to recovery, removal or death. The transition matrix \(V\) is given by
\[
V \;=\;
\begin{pmatrix}
\sigma + \mu & 0                & 0             & 0\\[4pt]
-\,\nu\sigma & \mu+d+\gamma     & 0             & 0\\[4pt]
(\nu - 1)\sigma & 0              & \mu+\delta+d    & 0\\[4pt]
0            & -\omega_s        & -\omega_a     & \varepsilon
\end{pmatrix}.
\]
Rows correspond to the ordering \((E,I_s,I_a,B)\). By construction the diagonal entries are the rates leaving each infected class.

\medskip
The next generation matrix $\psi_{NGM} = FV^{-1}$ and the basic reproduction number is the spectral radius of the next–generation matrix. Thus,
\[
\mathcal{R}_0 \;=\; \rho(\psi_{NGM})\;=\; \rho(FV^{-1}).
\]
We have 
\[
V^{-1} =
\begin{pmatrix}
\dfrac{1}{\sigma+\mu} & 0 & 0 & 0\\[10pt]
\dfrac{\nu\sigma}{(\sigma+\mu)(\mu+d+\gamma)} & \dfrac{1}{\mu+d+\gamma} & 0 & 0\\[10pt]
\dfrac{(1-\nu)\sigma}{(\sigma+\mu)(\mu+\delta+d)} & 0 & \dfrac{1}{\mu+\delta+d} & 0\\[10pt]
\dfrac{\nu\sigma\,\omega_s}{\varepsilon(\sigma+\mu)(\mu+d+\gamma)} 
+ \dfrac{(1-\nu)\sigma\,\omega_a}{\varepsilon(\sigma+\mu)(\mu+\delta+d)} &
\dfrac{\omega_s}{\varepsilon(\mu+d+\gamma)} &
\dfrac{\omega_a}{\varepsilon(\mu+\delta+d)} &
\dfrac{1}{\varepsilon}
\end{pmatrix}.
\]
Thus, we get the basic reproduction number 
\begin{equation}
\label{eqn:R0}
\mathcal{R}_0 = \frac{\sigma}{\sigma+\mu}\!\left[
\frac{\nu\,\beta_s}{\mu+d+\gamma} 
+ \frac{(1-\nu)\,\beta_a}{\mu+\delta+d}
+ \frac{\beta_B}{K\,\varepsilon}\!\left(
\frac{\nu\,\omega_s}{\mu+d+\gamma}
+ \frac{(1-\nu)\,\omega_a}{\mu+\delta+d}
\right)\!
\right].
\end{equation}

The prefactor \(\sigma/(\sigma+\mu)\) accounts for the probability that an exposed host survives the latent period and becomes infectious. The first two terms in the bracket are the sum of direct (host-to-host) contributions from symptomatic and asymptomatic infectives. The last two terms in the bracket quantifies the indirect (environmental) transmission that an infected host sheds virus into the environment at rates \(\omega_s,\omega_a\)), this virus decays at rate \(\varepsilon\), and environmental virus produces new infections at rate proportional to \(\beta_B\, S^*/(K\varepsilon)\).
The stochastic diffusion terms in the model do not change the location of the DFE (they vanish when \(E=I_s=I_a=B=0\)), but they do affect stability and sample-path behaviour near the DFE. The stochastic stability analysis near the equilibrium is discussed in next section.

\subsubsection{Endemic equilibrium (EE)}

We find the endemic equilibrium of the deterministic part of the system.
Set all stochastic terms and time derivatives to zero. When we analyse endemic equilibrium in epidemic models, we mean a steady state of the system where the number of individuals in each compartment remains constant over time. But, Brownian motion has no steady state: its increments are unbounded, so the noise term cannot vanish at a fixed positive level. The equilibrium \((S^{**},E^{**},I_s^{**},I_a^{**},R^{**},B^{**})\)
satisfies the algebraic system
\begin{align}
	0 &= \Lambda - \beta_a\frac{S^{**} I_a^{**}}{N^{**}} - \beta_s\frac{S^{**} I_s^{**}}{N^{**}} - \beta_B\frac{S^{**} B^{**}}{K+B^{**}} - \mu S^{**}, \label{eq:eq1}\\
	0 &= \beta_a\frac{S^{**} I_a^{**}}{N^{**}} + \beta_s\frac{S^{**} I_s^{**}}{N^{**}} + \beta_B\frac{S^{**} B^{**}}{K+B^{**}} - (\sigma+\mu) E^{**}, \label{eq:eq2}\\
	0 &= \nu\sigma E^{**} - (\mu+d+\gamma) I_s^{**}, \label{eq:eq3}\\
	0 &= (1-\nu)\sigma E^{**} - (\delta+\mu+d) I_a^{**}, \label{eq:eq4}\\
	0 &= \gamma I_s^{**} + \delta I_a^{**} - \mu R^{**}, \label{eq:eq5}\\
	0 &= \omega_s I_s^{**} + \omega_a I_a^{**} - \varepsilon B^{**}. \label{eq:eq6}
\end{align}
Here \(N^{**} = S^{**} + E^{**} + I_s^{**} + I_a^{**} + R^{**}\). We are supposed to find solution with \(E^{**}>0\)
(and hence \(I_s^{**},I_a^{**},B^{**},R^{**} \ge 0\)).\\
First solve the linear equations \eqref{eq:eq3}--\eqref{eq:eq6} in terms of \(E^{**}\). From equation \eqref{eq:eq3} and equation \eqref{eq:eq4},
\begin{align}
	I_s^{**} &= \alpha_s E^{**}, \label{eq:Is}\\
	I_a^{**} &= \alpha_a E^{**}. \label{eq:Ia}
\end{align}
where 
\[
\alpha_s \;:=\; \frac{\nu\sigma}{\mu+d+\gamma}, \qquad
\alpha_a \;:=\; \frac{(1-\nu)\sigma}{\delta+\mu+d}.
\]
are the two positive proportionality constants.\\
From equation \eqref{eq:eq5} we obtain
\[
R^{**} \;=\; \frac{\gamma I_s^{**} + \delta I_a^{**}}{\mu}
\;=\; \frac{\gamma \alpha_s + \delta \alpha_a}{\mu}\; E^{**}.
\]
Define for convenience
\[
\rho \;:=\; \frac{\gamma \alpha_s + \delta \alpha_a}{\mu},
\]
so \(R^{**} = \rho E^{**}\).
From equation \eqref{eq:eq6} we obtain
\[
B^{**} \;=\; \frac{\omega_s I_s^{**} + \omega_a I_a^{**}}{\varepsilon}
\;=\; \frac{\omega_s \alpha_s + \omega_a \alpha_a}{\varepsilon}\; E^{**}.
\]
Define
\[
\zeta \;:=\; \frac{\omega_s \alpha_s + \omega_a \alpha_a}{\varepsilon},
\]
so \(B^{**} = \zeta E^{**}\).
\vspace{6pt}
Now using the relations above we can express \(N^{**}\) in terms of \(S^{**}\) and \(E^{**}\). 
\[
N^{**} \;=\; S^{**} + E^{**} + \alpha_s E^{**} + \alpha_a E^{**} + \rho E^{**}\
\;=\; S^{**} + c\,E^{**}.
\]
where
\[
c \;:=\; 1 + \alpha_s + \alpha_a + \rho
\;=\; 1 + \alpha_s + \alpha_a + \frac{\gamma \alpha_s + \delta \alpha_a}{\mu}.
\]
Let us introduce two constants
\[
A_1 \;:=\; \beta_s \alpha_s + \beta_a \alpha_a,
\qquad
A_2 \;:=\; \beta_B \zeta .
\]
Here \(A_1\) collects the contributions of cattle-to-cattle transmission via symptomatic and asymptomatic infectives, and \(A_2\) collects the environmental contribution.\\
Now use the equation \eqref{eq:eq2} to relate \(S^{**}\) and \(E^{**}\).
Plugging \(I_s^{**}=\alpha_s E^{**}, I_a^{**}=\alpha_a E^{**}, B^{**}=\zeta E^{**}\) into \eqref{eq:eq2} gives
\[
\beta_a\frac{S^{**} \alpha_a E^{**}}{N^{**}} + \beta_s\frac{S^{**}\alpha_s E^{**}}{N^{**}}
+ \beta_B\frac{S^{**} \zeta E^{**}}{K+\zeta E^{**}}
\;=\; (\sigma+\mu) E^{**}.
\]
Factor \(E^{**}\) (and cancel \(E^{**}>0\)):
\[
S^{**} \left( \frac{A_1}{N^{**}} + \frac{A_2}{K+\zeta E^{**}}\right) \;=\; \sigma+\mu.
\]
Equivalently,
\[
\frac{\sigma+\mu}{S^{**}} \;=\; \frac{A_1}{N^{**}} + \frac{A_2}{K+\zeta E^{**}}.
\tag{15}
\]
This is one relation between \(S^*\) and \(E^*\) (recall \(N^*=S^*+cE^*\)).
Now, use the equation \eqref{eq:eq1} to express \(S^{**}\) by the infection pressure.
Let us define the instantaneous infection pressure (per susceptible)
\[
\lambda^{**} \;:=\; \frac{\beta_s I_s^{**}}{N^{**}} + \frac{\beta_a I_a^{**}}{N^{**}} + \frac{\beta_B B^{**}}{K+B^{**}}.
\]
Then equation \eqref{eq:eq1} is
\[
0 = \Lambda - S^{**}(\mu + \lambda^{**}),
\]
hence
\[
S^{**} \;=\; \frac{\Lambda}{\mu + \lambda^{**}}.
\tag{16}\label{eq:16}
\]
From equation \eqref{eq:eq2} we have \(\lambda^{**} S^{**} = (\sigma+\mu) E^{**}\) (because the LHS of equation \eqref{eq:eq2} equals \(\lambda^{**} S^{**}\));
therefore
\[
\lambda^{**} \;=\; \frac{(\sigma+\mu) E^{**}}{S^{**}}.
\tag{17}\label{eq:17}
\]
Substituting equation \eqref{eq:16} into equation \eqref{eq:17} we eliminate \(S^{**}\) and obtain an explicit expression for \(\lambda^{**}\)
as a function of \(E^{**}\):
\[
\lambda^{**}
= \frac{(\sigma+\mu) E^{**}}{\Lambda/(\mu+\lambda^{**})}
\;\Longrightarrow\;
\lambda^{**}(\Lambda - (\sigma+\mu)E^{**}) = (\sigma+\mu)\mu E^{**}.
\]
Hence, provided \(\Lambda - (\sigma+\mu)E^{**} \neq 0\),
\[
 \lambda^{**} \;=\; \frac{(\sigma+\mu)\mu\, E^{**}}{\Lambda - (\sigma+\mu)E^{**}}.
\tag{18}\label{eq:18}
\]
Above equation \eqref{eq:18} cancels the trivial solution \(E^{**}=0\); the endemic branch requires \(E^{**}>0\) and the denominator positive,
so a necessary admissibility restriction is \\ \(0<E^{**}<\dfrac{\Lambda}{\sigma+\mu}\).

\vspace{6pt}
Now obtain a single scalar equation for \(E^{**}\).

We have two expressions for the infection pressure \(\lambda^{**}\):\\
{(i)} By definition and using the proportionalities,
	\[
	\lambda^{**} \;=\; E^{**} \Big( \frac{A_1}{N^{**}} + \frac{A_2}{K+\zeta E^{**}} \Big),
	\tag{19}\label{eq:19}
	\]
	because \(I_s^{**}=\alpha_s E^{**},\,I_a^{**}=\alpha_a E^{**},\,B^{**}=\zeta E^{**}\) and \(N^{**}=S^{**}+cE^{**}\).\\
{(ii)} By elimination we have equation \eqref{eq:18} which gives \(\lambda^{**}\) solely as a function of \(E^{**}\). Equate the equations \eqref{eq:18} and \eqref{eq:19}; cancel \(E^{**}\) (since we search for \(E^{**}>0\)):
\[
\frac{(\sigma+\mu)\mu}{\Lambda - (\sigma+\mu)E^{**}}
\;=\; \frac{A_1}{N^{**}} + \frac{A_2}{K+\zeta E^{**}}.
\tag{20}\label{eq:20}
\]
We have \(N^{**} = S^{**} + cE^{**}\) and \(S^{**} = \Lambda/(\mu+\lambda^{**})\) with \(\lambda^{**}\) given by equation \eqref{eq:18}.
Thus the left-hand side is an explicit function of \(E^{**}\), and the right-hand side is also an explicit function of \(E^{**}\) after substitution.
Equation \eqref{eq:20} is the required scalar (nonlinear) equation for \(E^{**}\).

Now \(E^{**}\) is known, we can recover all remaining equilibrium components explicitly. Suppose \(E^{**}>0\) is a (positive) solution of equation \eqref{eq:18} that satisfies \(E^{**}<\Lambda/(\sigma+\mu)\). Then compute
\[
\lambda^{**}=\; \frac{(\sigma+\mu)\mu\, E^{**}}{\Lambda - (\sigma+\mu)E^{**}}\quad\text{from equation } \eqref{eq:20},
\qquad S^{**}=\frac{\Lambda}{\mu+\lambda^{**}},
\qquad N^{**}=S^{**}+cE^{**}.
\]
Then the remaining components are
\begin{align*}
	I_s^{**} &= \alpha_s E^{**},\\
	I_a^{**} &= \alpha_a E^{**},\\
	R^{**}   &= \rho E^{**} = \frac{\gamma \alpha_s + \delta \alpha_a}{\mu}\,E^{**},\\
	B^{**}   &= \zeta E^{**} = \frac{\omega_s \alpha_s + \omega_a \alpha_a}{\varepsilon}\,E^{**}.
\end{align*}
All of these are nonnegative provided \(E^{**}\ge0\). The admissibility condition for this analysis is \(0<E^{**}<\dfrac{\Lambda}{\sigma+\mu}\),
which ensures the denominator in equation \eqref{eq:18} is positive and \(S^{**}>0\).

\section{Numerical simulations} \label{sec:numel}
Here, we illustrate the theoretical results from the preceding section through numerical simulations of models \eqref{eq:ode-full} and \eqref{eq:sde-full}. Unless specified otherwise, all simulations adopt the parameter values in Table \ref{tab:parameterdescr}. For numerical simulation of the stochastic AI model, we employed the Euler--Maruyama scheme, which is the standard first–order discretization method for It\^o stochastic differential equations \cite{Kloeden2012}. At each time step, the deterministic drift terms of the model were updated using $\Delta t$, while the stochastic perturbations were incorporated through multiplicative noise components of the form $\sigma_X X_t\,\Delta W_t$. 
Let $\Delta t > 0$ be the step size and define discrete times
\[
t_k = k\Delta t,\qquad k=0,1,\ldots,N,
\]
where $N = T/\Delta t$ is the total number of steps up to a final simulation time $T$ (500 days).
For a Wiener increment we use the discrete approximation
\[
\Delta W_k = W(t_{k+1}) - W(t_k) 
\sim \mathcal{N}(0,\Delta t),
\]
and in computations we have
\[
\Delta W_k = \sqrt{\Delta t}\, Z_k,\qquad Z_k\sim\mathcal{N}(0,1).
\]
Let $X_k$ denote the numerical approximation of $X(t_k)$. Euler--Maruyama update Rule applied to equation (\ref{eq:Ito SDE}) is
\begin{equation*}
	X_{k+1} 
	= 
	X_k 
	+ f_X(X_k)\,\Delta t 
	+ \sigma_X X_k\, \Delta W_k.
	\label{eq:EM-update}
\end{equation*}

\begin{figure*}[h!]
\centering
 \includegraphics[width=\textwidth]{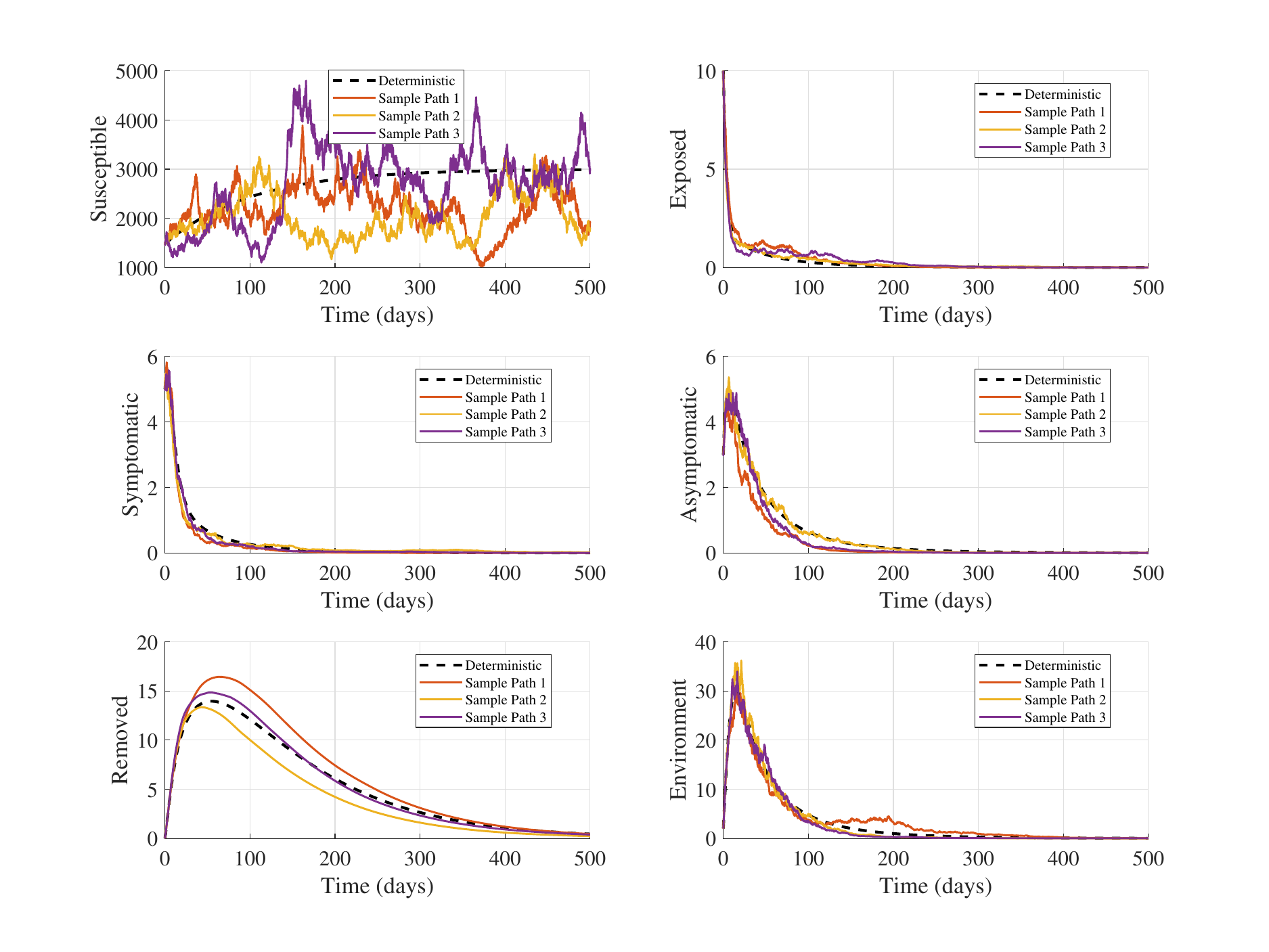}
\caption{Deterministic and stochastic trajectories of model \eqref{eq:sde-full} under disease-free equilibrium conditions ($\mathcal{R}_0 < 1$). All parameter values used in this simulation are listed in Table \ref{tab:parameters}.}
\label{fig:DFE}
 \end{figure*}

The numerical simulation for the DFE analysis shown in Figure \ref{fig:DFE}, compares the deterministic solution of the model \eqref{eq:sde-full} with a stochastic realization that incorporates multiplicative environmental and demographic noise. Under the chosen parameter set (Table~\ref{tab:parameters}) corresponding to $\mathcal{R}_0 < 1$, the system approaches to disease-free equilibrium, and all infected compartments vanish over time.
The deterministic trajectories display smooth exponential decay in the exposed and infectious classes and a monotone return of the susceptible population to its disease-free level. The stochastic realization follows the same qualitative pattern but exhibits sample-path fluctuations. While early dynamics show small irregular oscillations around the deterministic decline in $E, I_a, I_s$ and $B$, these perturbations do not lead to secondary outbreaks.

\begin{table}[H]
\centering
\footnotesize
\caption{Parameters of the stochastic AVI model \ref{eq:stochastic-system}, their values, and sources.}
\label{tab:parameterdescr}

\begin{tabular}{ccc|ccc|ccc}
\toprule
\textbf{Parameter} & \textbf{Value} & \textbf{Source} &
\textbf{Parameter} & \textbf{Value} & \textbf{Source} &
\textbf{Parameter} & \textbf{Value} & \textbf{Source} \\
\midrule

$\Lambda$      & 30        & Assumed
& $K$           & 500      & Assumed
& $\gamma$      & 0.10     & \cite{Bellotti2024,Rawson2025} \\

$\beta_a$      & 0.004  &  Assumed
& $\beta_s$     & 0.005     & Assumed
& $\beta_B$     & 0.002  & \cite{Bellotti2024,Rawson2025} \\

$\sigma$       & 0.20     & \cite{Rawson2025}
& $\nu$         & 0.50     & Assumed
& $\delta$      & 0.05     & Assumed \\

$d$            & 0.01     & \cite{Regassa2024}
& $\omega_s$    & 0.5     & Assumed
& $\omega_a$    & 0.4     & Assumed \\

$\epsilon$     & 0.10     & Assumed
& $\mu$         & 0.01     & \cite{WaFe08}
& $\sigma_s$    & 0.05     & Assumed \\

$\sigma_S$     & 0.05  & Assumed
& $\sigma_e$    & 0.05     & \cite{Rawson2025}
& $\sigma_B$    & 0.05     & Assumed \\

\bottomrule
\end{tabular}
\end{table}

Across all infected compartments $(E, I_s, I_a)$ and the environmental viral load $(B)$, both deterministic and stochastic paths rapidly decline toward zero, confirming the inability of the pathogen to persist within the host population and also in the environment. The small stochastic fluctuations do not enable infection rebound.
For susceptible population, despite large sample-path variability, the stochastic path does not show a sustained epidemic wave. It fluctuates around the DFE and remains bounded. Both infectious classes, symptomatic and asymptomatic, decay rapidly in the stochastic simulation, consistent with a negative noise-adjusted growth rate and confirming that random perturbations are insufficient to sustain transmission. Stochastic noise causes small transient departures from the deterministic decline, but the exponential decay dominates. A transient peak exists in the environment early on as infectious hosts shed pathogen into the environment, then decays toward zero. Environmental loading transiently rises but is insufficient to maintain transmission long term due to environmental decay and also low infection prevent sustained environmental forcing. Thus, when $(\mathcal{R}_0 < 1)$, neither host-to-host transmission nor environmental transmission can sustain endemic infection in cattle under the current parameter settings.

\begin{figure*}[h!]
\centering
 \includegraphics[width=\textwidth]{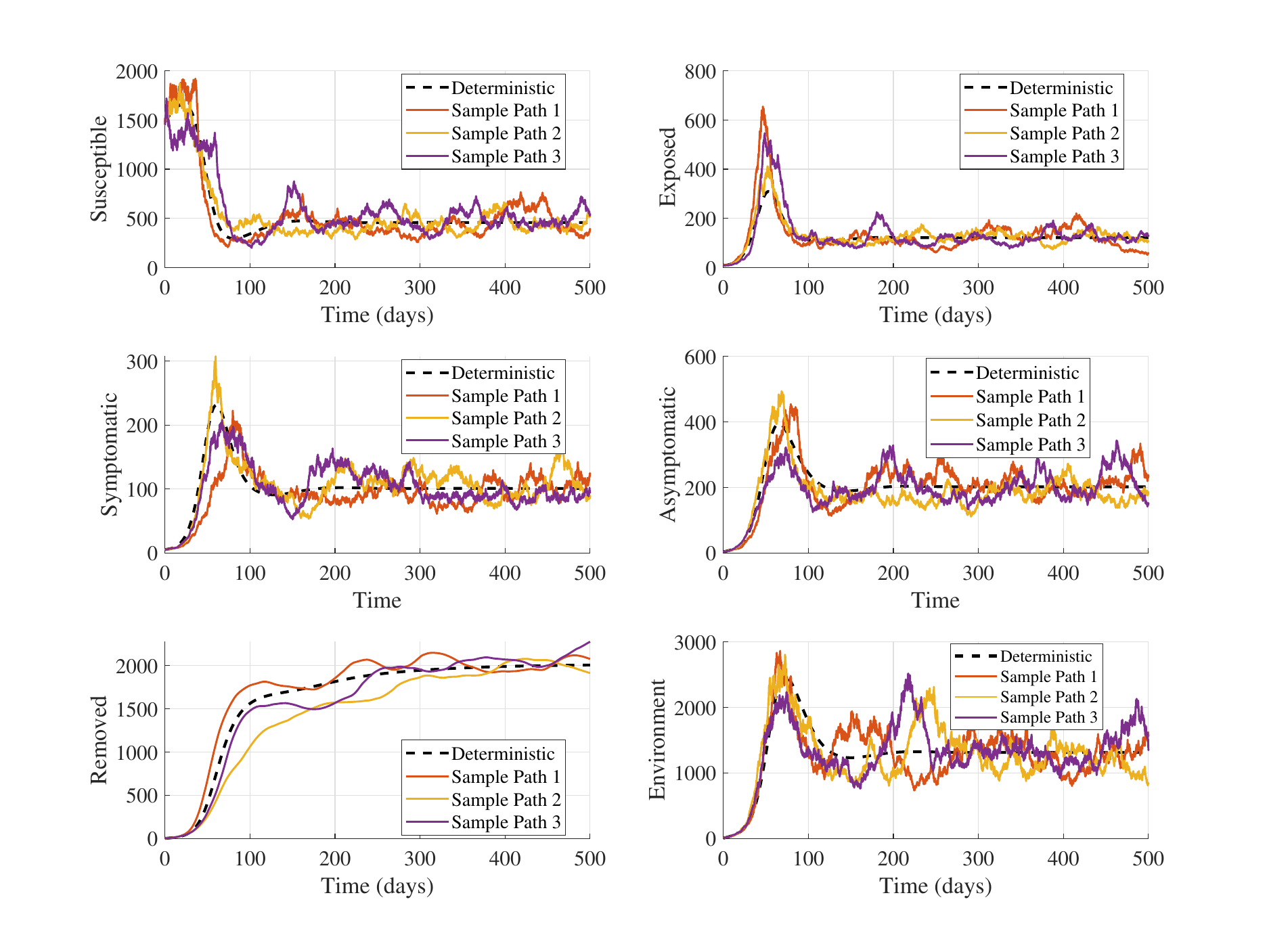}
\caption{Transmission dynamics of model (\ref{eq:sde-full}) under endemic equilibrium conditions ($\mathcal{R}_0 > 1$). The black dashed curve in each panel represents the deterministic trajectory while colored solid lines represent three stochastic realizations driven by independent Wiener processes.}
 \label{fig:EE}
 \end{figure*}

Figure \ref{fig:EE} illustrates how the model behaves under an endemic equilibrium over a 500-day period. All model parameters remain the same as in the baseline disease-free analysis, except for an increase in the asymptomatic transmission rate $\beta_a$, which elevates the basic reproduction number to $\mathcal{R}_0 = 3.1945$, pushing the system from a disease-free regime into a stable endemic state.
This modification strengthens the indirect and silent spread of infection within the host population, resulting in a higher basic reproduction number. Now the disease no longer dies out; instead, it persists in the population, and the system settles into an endemic state.

In this model, asymptomatic infectious individuals play a crucial role because they continue to interact with susceptible hosts without displaying clinical symptoms. As a result, they remain undetected, untreated, and unrestricted in movement, creating a hidden but highly effective transmission pathway. When $\beta_a$ is increased, these individuals contribute more strongly to new infections, amplifying the force of infection. The results shown in Figure \ref{fig:EE} also support the condition for \(0<E^{**}<\dfrac{\Lambda}{\sigma+\mu}\) in the analysis of endemic equilibrium based on the parameter values as defined in Table \ref{tab:parameterdescr}.

\subsection{Sensitivity analysis}
\begin{figure}[h!]
    \centering
    \includegraphics[width=0.7\linewidth]{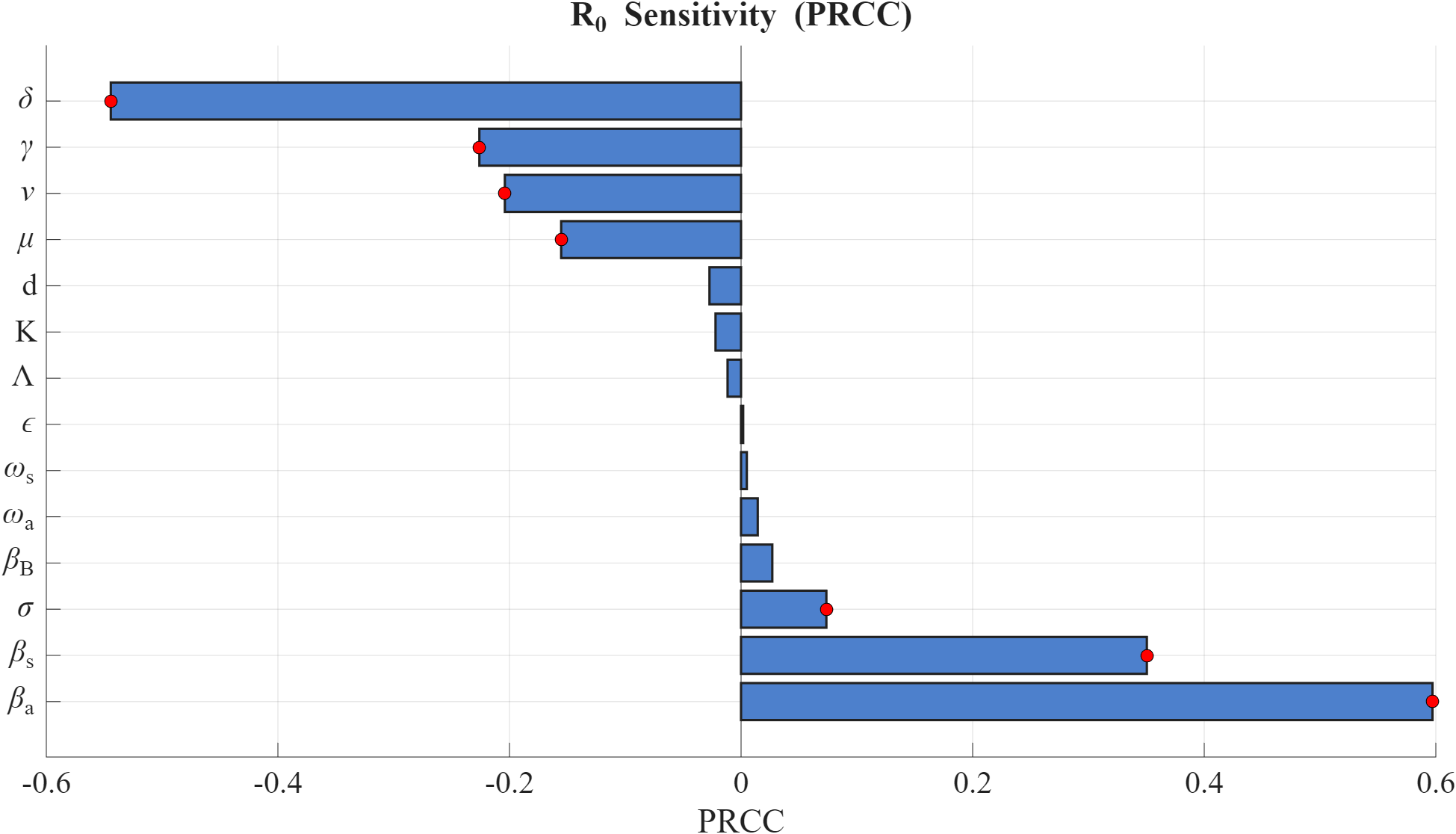}
    \caption{Global sensitivity analysis of $\mathcal{R}_{0}$ with respect
to model \eqref{eq:sde-full} parameters using Partial Rank Correlation Coefficients (PRCC) and Latin Hypercube Sampling (LHS).}
    \label{fig:sensitivity}
\end{figure}

\begin{figure*}[h!]
\centering
   \begin{subfigure}[b]{.4\linewidth}
   \centering
   \caption{}
\includegraphics[width=\textwidth]{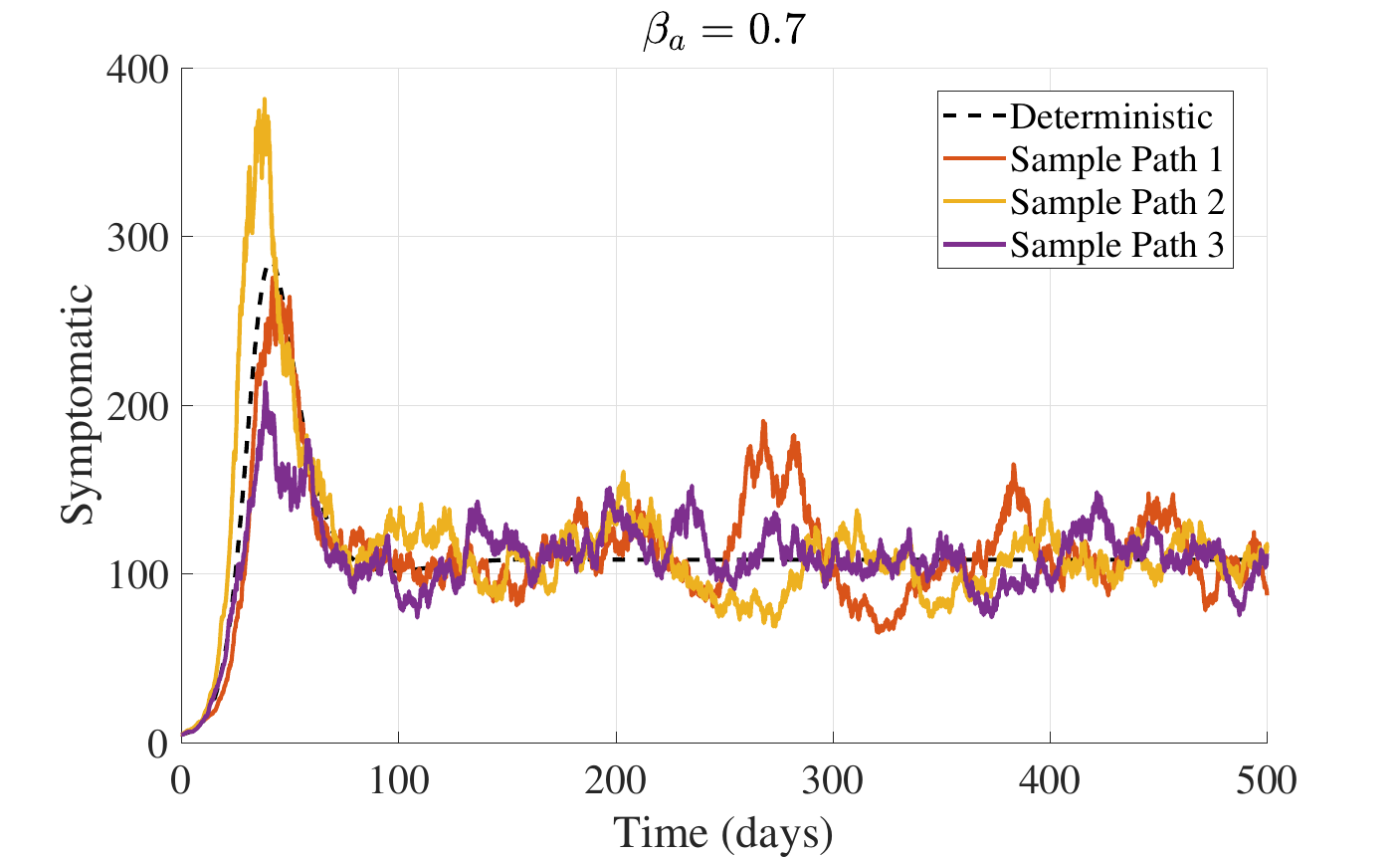}
     \label{fig:beta07}
  \end{subfigure}%
 \begin{subfigure}[b]{.4\linewidth}
   \centering
   \caption{}
\includegraphics[width=\textwidth]{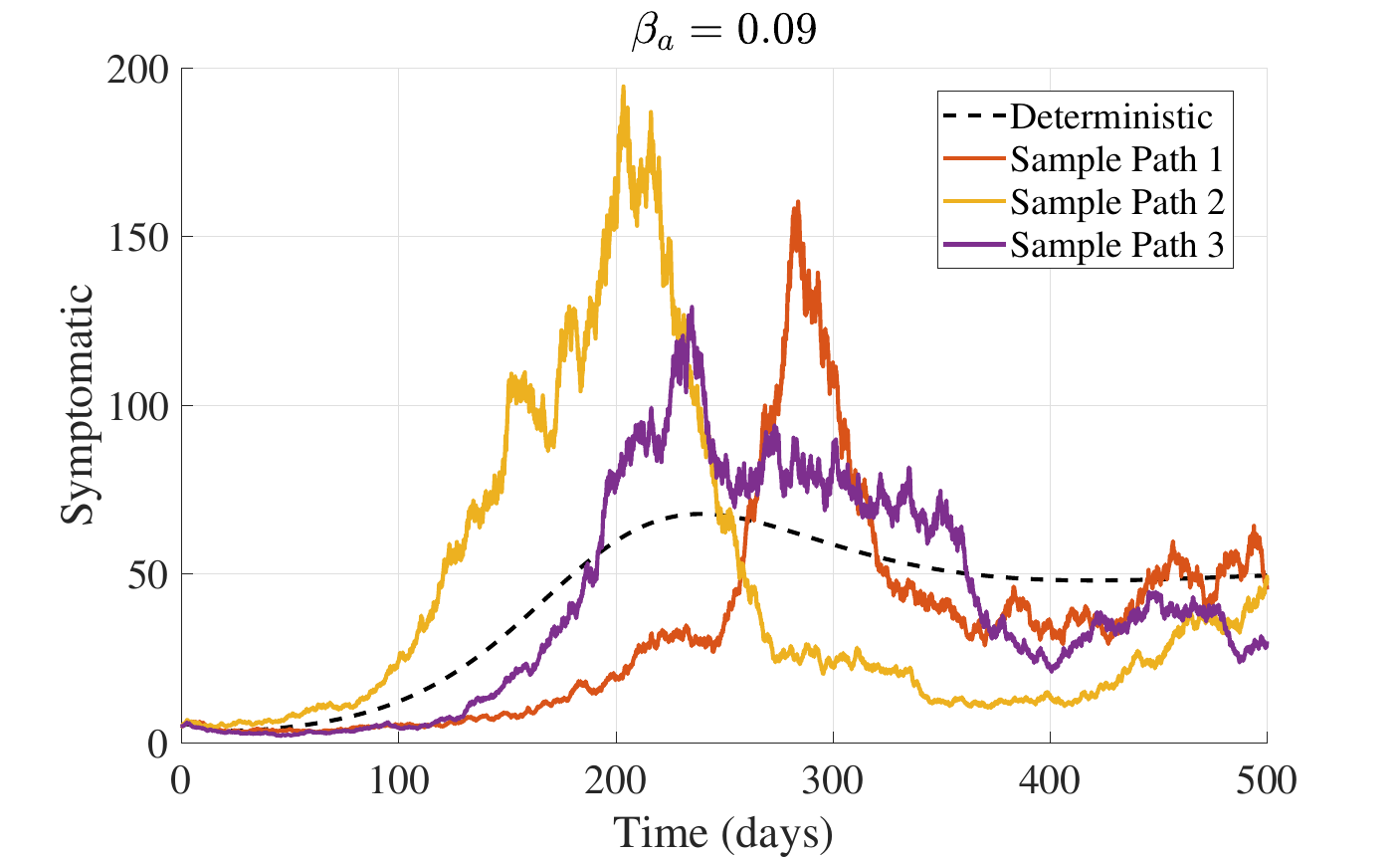}
     \label{fig:betaa009}
  \end{subfigure}%
   \caption{Epidemic trajectories of model \eqref{eq:sde-full} are presented for $\beta_a = 0.7$ and $\beta_a = 0.09$, with all other parameters taken as in Table~\ref{tab:parameterdescr}. The deterministic solution is depicted by dashed lines, while while the stochastic sample paths are shown as solid lines.}
    \label{fig:twobetaa}
  \end{figure*}

\begin{figure*}[h!]
\centering 
\begin{subfigure}[b]{.4\linewidth}
   \centering
   \caption{}
\includegraphics[width=\textwidth]{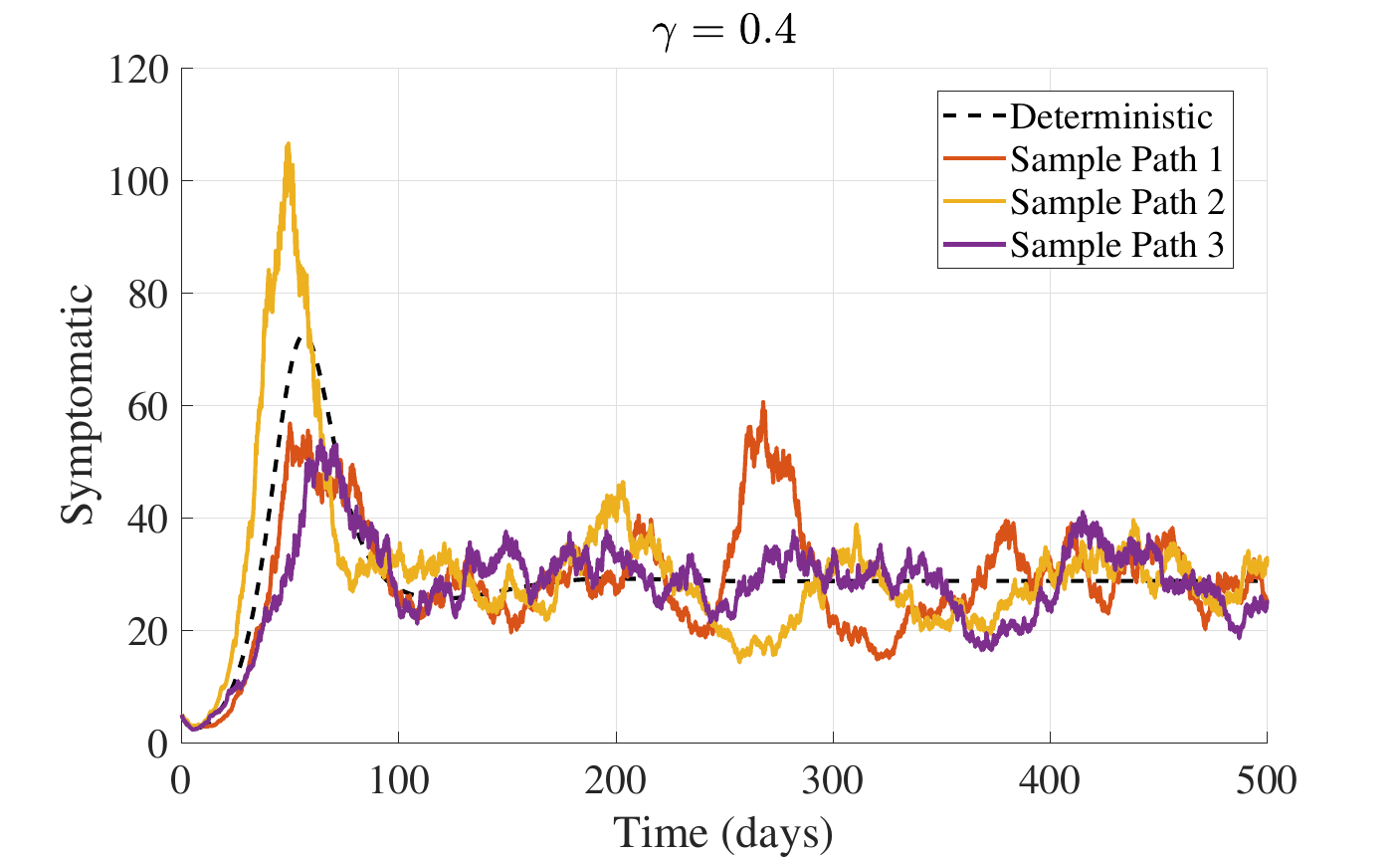}
     \label{fig:gamma04}
  \end{subfigure}%
 \begin{subfigure}[b]{.4\linewidth}
   \centering
   \caption{}
\includegraphics[width=\textwidth]{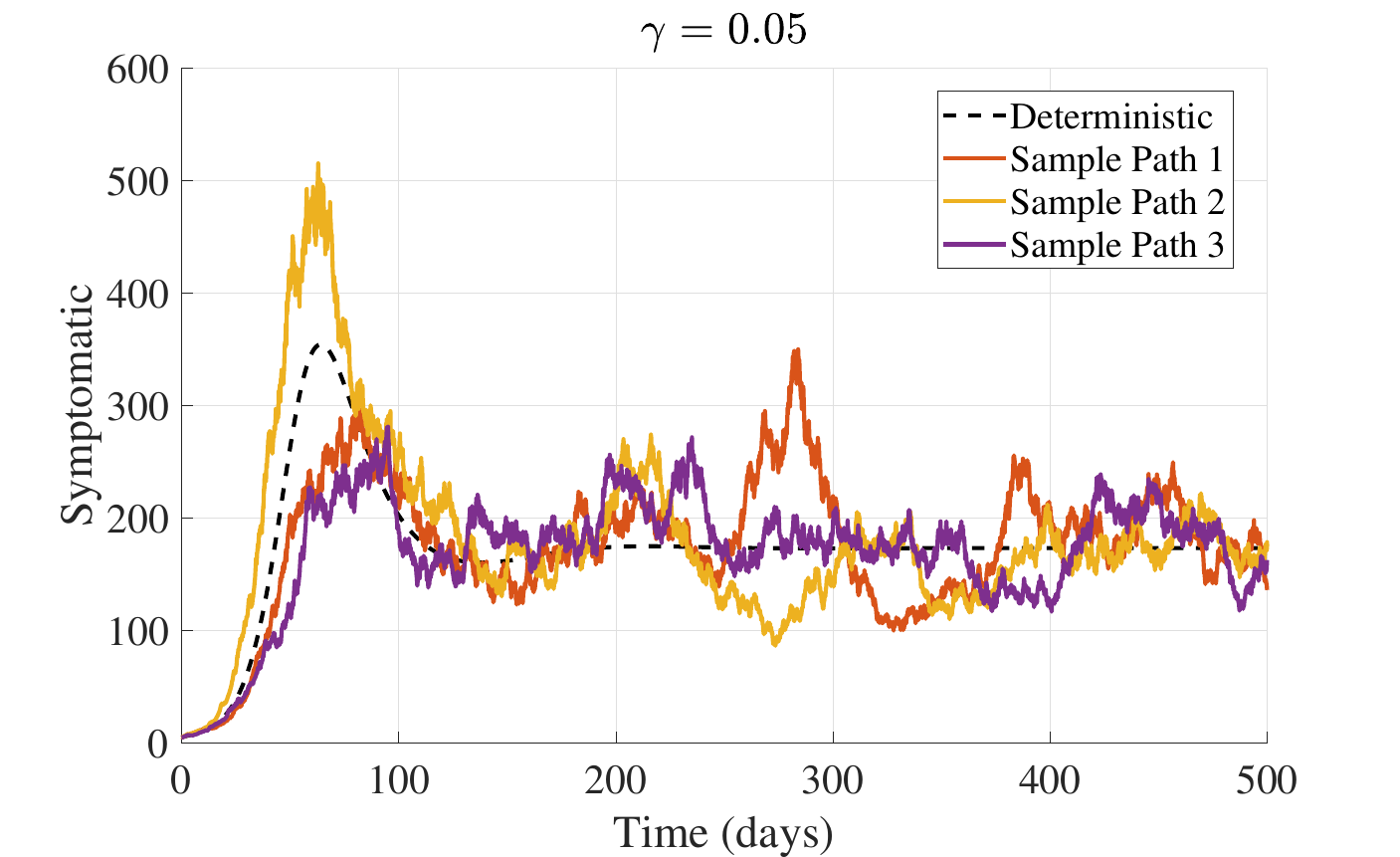}
     \label{fig:gamma005}
  \end{subfigure}%
   
   \caption{Epidemic trajectories of model \eqref{eq:sde-full} are presented for $\gamma=0.4$ and $\gamma=0.05$, with all other parameters taken as in Table~\ref{tab:parameterdescri}. The deterministic solution is depicted by dashed lines, while the stochastic sample paths are shown as solid lines.}
    \label{fig:twogamma}
  \end{figure*}

The extinction or persistence of the epidemic depends on the threshold quantity $\mathcal{R}_{0}$. We performed a global sensitivity analysis of $\mathcal{R}_{0}$, as defined in \eqref{eqn:R0}, using Latin Hypercube Sampling (LHS) and Partial Rank Correlation Coefficient (PRCC) techniques \cite{Abidemi22, Reju24}. The results are presented in Figure \ref{fig:sensitivity}, where the bars marked with a red dot indicate the parameters that exert a significant influence on $\mathcal{R}_{0}$. A positive (negative) PRCC value indicates a direct (inverse) relationship with $\mathcal{R}_{0}$, meaning that an increase in the parameter results in an increase (decrease) in $\mathcal{R}_{0}$. For instance, $\beta_{a}, \beta_{s}$, and $\sigma$ have positive values, so increasing any of these parameters increases the value of $\mathcal{R}_{0}$. This suggests that larger values of these parameters contribute to a higher infection burden. Conversely, the parameters $\delta,~\gamma$, and $\nu$ have negative PRCC values, thus, increasing their values results in reduction of $\mathcal{R}_{0}$. Some results of the sensitivity analysis are in simulations in Figures~\ref{fig:twobetaa} and \ref{fig:twogamma} for two values of $\beta_a$ and $\gamma$, respectively. 
Figure~\ref{fig:twobetaa} presents model simulations for two different values of asymptomatic transmission rate $\beta_a$. When $\beta_a = 0.7$, the deterministic trajectory reaches a peak of about 280 symptomatic cattle within roughly 70 days, and the stochastic realizations follow a similarly rapid rise. In contrast, for $\beta_a = 0.09$, the deterministic peak is lower, at around 65 symptomatic cattle, and occurs later, between 200 and 250 days, with the stochastic trajectories exhibiting greater variability. Consequently, higher transmission leads to earlier and larger outbreaks, whereas lower transmission results in smaller, delayed peaks accompanied by increased stochastic variability.

As illustrated in Figure~\ref{fig:twogamma}, changing the recovery rate $\gamma$ markedly alters the course of the epidemic. For $\gamma = 0.4$, the deterministic trajectory peaks at about day 72 with roughly 75 symptomatic cattle, and the stochastic sample paths remain closely grouped around this relatively small outbreak, indicating rapid infection clearance. In contrast, for $\gamma = 0.05$, the deterministic curve increases more gradually but attains a much higher peak of approximately 350 symptomatic cattle around day 85, with the corresponding sample paths displaying greater dispersion. Consequently, a higher recovery rate reduces epidemic magnitude and shortens its duration, whereas a lower recovery rate extends the infectious period, intensifies transmission, and yields larger, later peaks.

\subsection{Varying noise intensity}

\begin{figure*}[h!]
\centering
 \begin{subfigure}[b]{.4\linewidth}
   \centering
   \caption{}
\includegraphics[width=\textwidth]{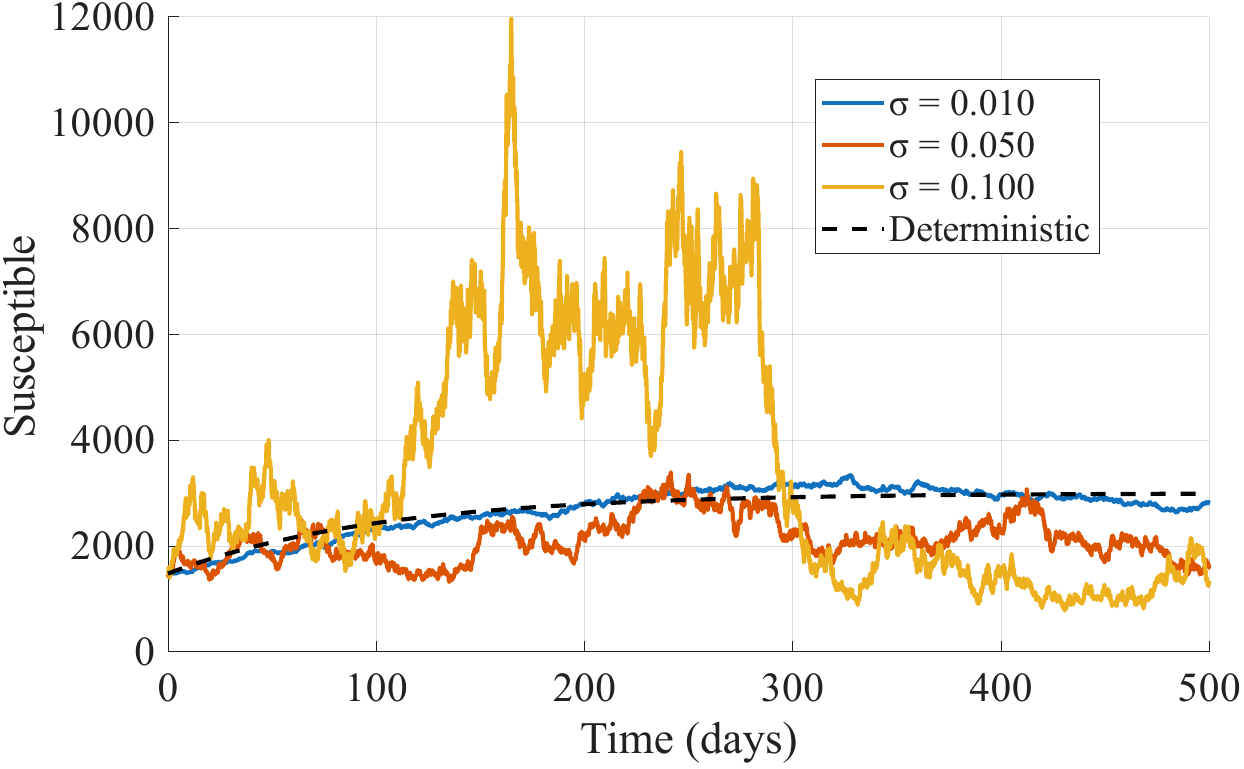}
     \label{fig:S_VarySig}
  \end{subfigure}\\%
 \begin{subfigure}[b]{.4\linewidth}
   \centering
   \caption{}
\includegraphics[width=\textwidth]{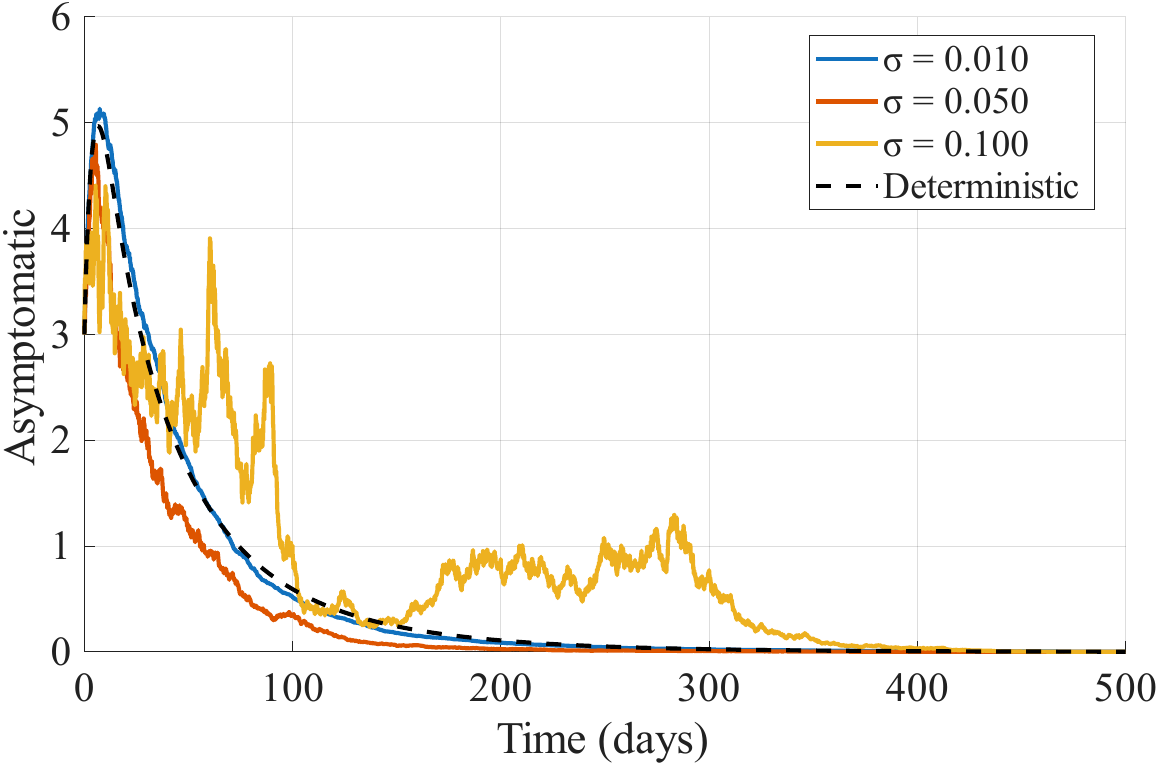}
     \label{fig:Ia_varySigma}
  \end{subfigure}%
  \begin{subfigure}[b]{.4\linewidth}
\centering
\caption{}
\includegraphics[width=\textwidth]{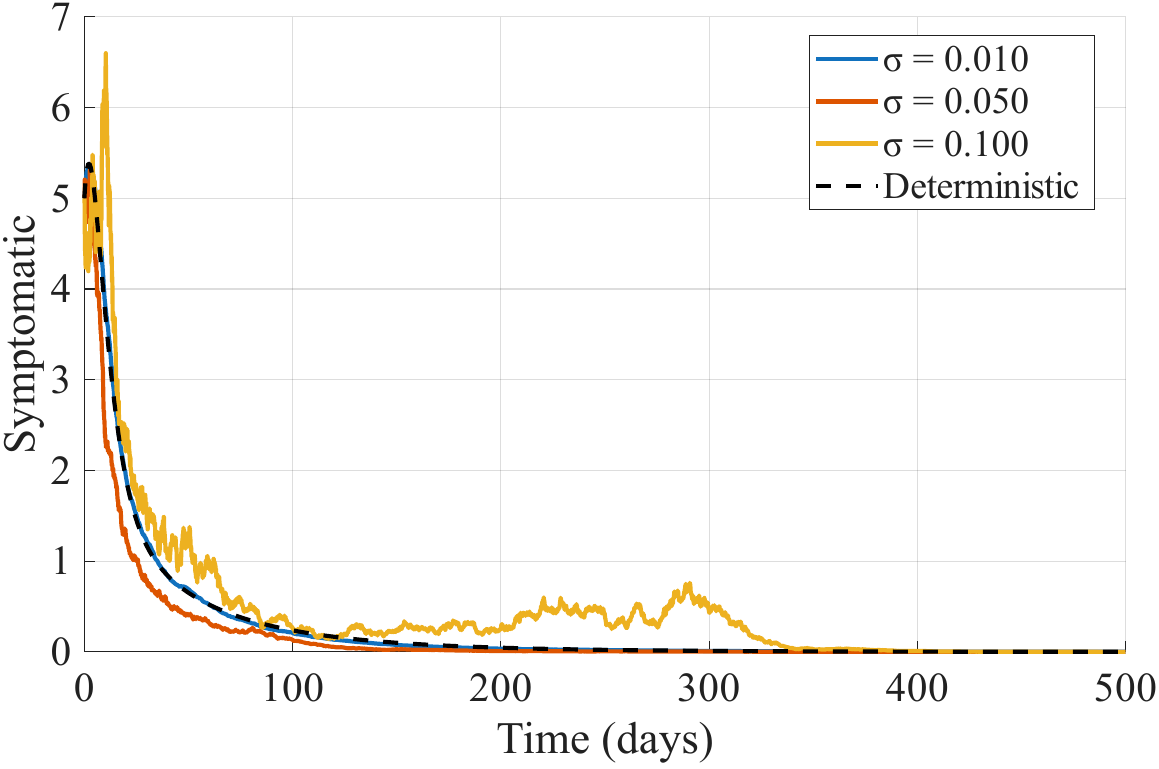}
\label{fig:Is_varySigma}
 \end{subfigure}%
   \caption{Epidemic trajectories under varying noise intensity, $\sigma =0.010,0.050,$ and $0.100$, with all other parameters taken as in Table~\ref{tab:parameterdescr}. The deterministic solution is depicted by dashed lines, while the stochastic sample paths are shown as solid lines.}
    \label{fig:varySigma}
  \end{figure*}

Here we examine how changes in the noise intensity, $\sigma_k$ (for $k=S,E,I_a,I_s$), affect the epidemic dynamics. The outcomes for three distinct noise intensities are presented in Figure~\ref{fig:varySigma}. When noise is low ($\sigma_k = 0.010$), the system behaves almost deterministically: the epidemic trajectories for the asymptomatic, symptomatic, and susceptible compartments are smooth, with clearly defined and predictable peaks. The final epidemic size remains tightly clustered around the deterministic value, and the timing of the main outbreak is highly reproducible across realizations. As the noise level increases to $\sigma_k = 0.050$, moderate stochastic fluctuations introduce visible variability among sample paths. The peak incidence in both asymptomatic and symptomatic groups may be either enhanced or attenuated compared with the deterministic scenario, and the timing of the epidemic peak becomes more variable. Consequently, the distribution of outbreak sizes and durations broadens, complicating both epidemic forecasting and the timely deployment of control strategies.

At high noise levels ($\sigma_k = 0.100$), stochasticity strongly shapes the epidemic course. Sample paths can exhibit multiple irregular peaks, extended low-level circulation, or even premature fade-outs. The symptomatic compartment, in particular, may show large, erratic surges that are not captured by the deterministic model. The susceptible population also displays a much more uncertain trajectory, with some realizations experiencing rapid depletion and others retaining a substantial susceptible pool for long periods. Overall, increasing noise intensity drives the system away from deterministic predictability and toward a regime where epidemic size, peak magnitude, and timing become highly variable and uncertain.

\section{Conclusion} \label{sec:conclusion}
This study developed and analysed a stochastic model to explore the transmission dynamics of avian influenza (AI) in cattle populations with both direct and environmental routes of infection. Although AI infections in dairy cattle have only recently gained international attention, especially after the 2024 H5N1  detections in North American herds \cite{Rawson2025}, this model captures the dual role of direct contact transmission and indirect exposure through contaminated environment. 

A key insight from the model is the contribution of asymptomatic infections. As seen in livestock influenza systems \cite{Mrope2024,sreenivasan2024emerging}, a substantial proportion of infections may be mild or subclinical, yet still infectious. This justifies the separation of symptomatic and asymptomatic infectious states. The model shows that the asymptomatic class can sustain transmission even when symptomatic prevalence is low, particularly when coupled with environmental shedding as seen in Figures \ref{fig:Ia_varySigma} and \ref{fig:Is_varySigma}. This aligns with recent reports suggesting subclinical AI infections may occur in cattle herds \cite{Mrope2024}, making detection and control more challenging. 

As seen in Figure~\ref{fig:EE}, an increase in the value of the asymptomatic transmission rate $\beta_a$ leads to the model showing sustained infection and a stable endemic state. The results agree with the sensitivity analysis results in Figure~\ref{fig:sensitivity}. This is biologically feasible for cattle herds, especially in situations with repeated contamination of shared water sources, milking systems, holding yards, or slurry pits.

The introduction of stochasticity into the model plays a crucial role in capturing the inherent variability present in real livestock systems. In contrast to purely deterministic dynamics, which predict smooth trajectories and fixed thresholds, the stochastic model reveals how random fluctuations can substantially alter both the short-term patterns and long-term outcomes of avian influenza transmission in cattle population (see Figure \ref{fig:twobetaa}). The simulation results show that inclusion of randomness in cattle compartments make the system fluctuates around the endemic equilibrium as seen in Figure \ref{fig:twogamma}. This is important because cattle farms may differ in herd size, movement, milking schedule, and hygiene, all of which introduce randomness that deterministic models cannot capture. 

Varying the noise intensities across the compartments generally amplifies variability in disease prevalence, but the effects are not uniform across compartments. Instead, the combined action of noise modifies how infection enters, circulates, and persists in the herd. When the overall noise level is low, the stochastic trajectories remain close to the deterministic solution: outbreaks follow predictable curves, and the infection either dies out or settles around a stable endemic level depending on the value of $\mathcal{R}_0$. However, as noise intensities increase, the system becomes progressively more irregular, and deviations from the deterministic baseline are amplified, as seen in Figure~\ref{fig:varySigma}.

One of the most important consequences of higher stochasticity is the increased likelihood of early extinction when the system is near the epidemic threshold. When $\mathcal{R}_0 > 1$ and with varying noise intensities, the infection does not approach a single equilibrium. Instead, the system fluctuates around the endemic level, sometimes producing large transient spikes driven by random increases in viral shedding or environmental contamination. This behaviour is biologically feasible because cattle operations are subject to changes in weather affecting virus survival, variability in shedding between cattle, irregular contact patterns, and differences in cleaning routines; all of these factors contribute to environmental and demographic noise \cite{Zhang2019PhysA}.

These results show the importance of stochastic effects when assessing outbreak risks in cattle settings. In real situations, herds rarely exhibit smooth, deterministic dynamics; instead, their disease patterns reflect the cumulative effect of multiple sources of randomness. Incorporating noise into the model results in a more realistic representation of AI dynamics and highlights the importance of measures that reduce environmental variability, improve hygiene, and limit contact patterns in case of an outbreak.

Our findings of with-in herd analysis suggested that transmission rates of asymptomatic and symptomatic cattle, and noise intensity played a significant role in the spread of HPAI disease. This research work is not exhaustive, but it is expected that this would be helpful for livestock farmers, health professionals, and policymakers to take necessary measures against the transmission dynamics of HPAI.

\section{Limitations and future scope}\label{sec:limitation}
Although the stochastic model developed in this study provides useful insights into the transmission dynamics of HPAI in cattle population, it has some limitations. The model simplifies heterogeneity between animals by assuming homogeneous mixing, and uniform susceptibility, exposure, and shedding rates. In practice, cattle may differ in age, behaviour, housing conditions, immune responses, and viral shedding patterns. Such heterogeneity can significantly affect transmission, particularly for asymptomatic carriers who may shed virus intermittently or at lower levels. Incorporating age-structured or behaviour-based subpopulations could therefore lead to more accurate predictions. Also, several parameter values were selected from literature on avian influenza in other species due to limited data on cattle-specific AI dynamics. While this approach provides biologically reasonable ranges, it inevitably introduces uncertainty. Direct parameterisation using cattle outbreak data would improve the reliability and predictive accuracy of the model, but such datasets are currently limited.


In future, we will extend this model by incorporating spatially structured systems to capture heterogeneity between farms. Also, we will consider integrating available HPAI longitudinal data to obtain more accurate and context-dependent parameter estimates. Another important extension is evaluating control strategies, such as vaccination and movement restrictions. Stochastic optimal control and cost-effectiveness frameworks could help in understanding avian influenza transmission in cattle and support more resilient disease management strategies.

\section*{Author contributions}

Conceptualization, P.T and H.O.F.; methodology, P.T., M.S. and H.O.F.; software, P.T., M.S. and H.O.F.; validation, P.T. and H.O.F.; formal analysis, P.T. and M.S.; investigation, P.T., M.S. and H.O.F.; resources, P.T., M.S. and H.O.F.; writing---original draft preparation, M.S. and P.T.; writing---review and editing, P.T. and H.O.F.; visualization, P.T., M.S. and H.O.F.; supervision, P.T. and H.O.F.; project administration, P.T. and H.O.F.; funding acquisition, P.T. and H.O.F.

\section*{Funding}
This research was funded by AUT School of Engineering, Computer and Mathematical Sciences Research Grant.

\section*{Code Availability}
The code used to generate the results and figures in this study is available at \url{https://github.com/hamfat/Stochastic-HPAI-Model}

\end{document}